\documentclass{ws-procs9x6}

\def\CH{\mathcal{H}}
\def\CD{\mathcal{D}}
\def\CC{\mathcal{C}}
\def\CM{\mathcal{M}}
\def\CO{\mathcal{O}}
\def\R{\mathbb{R}}
\def\P{\mathbb{P}}
\def\d{\partial}
\def\scal#1{\langle#1\rangle}
\let\eps\varepsilon
\def\Lip{\mathrm{Lip}}

\begin{document}

\title{Hypoellipticity in Infinite Dimensions}

\author{Martin Hairer}

\address{Courant Institute, NYU and University of Warwick}

\begin{abstract}
We consider semilinear parabolic stochastic PDEs driven by additive noise.
The question addressed in this note is that of the regularity of transition probabilities. If the
equation satisfies a H\"ormander `bracket condition', then any finite-dimensional
projection of the solution has a smooth density with respect to Lebesgue measure.
One key ingredient in the argument is a bound on `Wiener polynomials' that plays a
role analogue to Norris' lemma.
\end{abstract}

\section{Introduction}

In this note, we report on recent results obtained in collaboration with J.C. ~Mattingly\cite{Hypo}
regarding the behaviour of transition probabilities for a large class of semilinear stochastic PDEs. While these
results were motivated mainly by the study of the long-time behaviour of solutions\cite{HairerMattingly06AOM}, 
the aspect that we
will focus on in this note is that of the regularity of transition probabilities which is an interesting
(and still partially open) mathematical question in its own right. 

The class of problems that we are able to treat are semilinear parabolic stochastic PDEs of the type
\begin{equation}\label{e:main}
du = -Au\,dt + F(u)\,dt + Q\,dW(t)\;,\quad u_0 \in \CH\;,
\end{equation}
where $\CH$ is a separable Hilbert space, $A$ is a positive selfadjoint negative linear operator with compact resolvent,
$W$ is a finite-dimensional Wiener process taking values in $\R^d$ and $Q\colon \R^d \to \CH$ with $Qe_i = q_i$. 
The nonlinearity 
$F \colon \CD(A^\alpha) \to \CH$ is assumed to be  `polynomial' in the sense that it can be written as
a finite sum of multilinear terms. Prime examples of equations satisfying our assumptions are the 2D stochastic
Navier-Stokes equations on a bounded regular domain, reaction-diffusion equations in dimension
$3$ or less, the stochastic Kuramoto-Sivashinsky equation, stochastic Burgers' equation, etc.

We address the question whether the solutions to \eref{e:main} have `smooth' transition probabilities.
For finite-dimensional SDEs, the meaning of this question is clear: do the transition probabilities have a 
$\CC^\infty$ density with respect to Lebesgue measure? In the infinite-dimensional case, it is much less
clear what we mean by `smooth' since there is no natural reference measure. In the special case where
$W$ is a cylindrical Wiener process on $\CH$ and $Q$ is `large' (either with bounded inverse or
comparable to some inverse power of $A$), it is often possible\cite{ZDP,FM,b:GoldysMaslowski06} 
to show that the transition
probabilities have a density with respect to the Gaussian measure which is invariant for the linearised equation
(i.e. the same equation with $F=0$) and that this density is well-behaved. However, we are interested in the case
where $Q$ has finite-dimensional range, so that these arguments certainly do not apply. The next best notion of 
smoothness is then the following: 

\begin{question}\label{mainquest}
Given $N>0$ and an orthonormal projection $\Pi \colon \CH \to \R^N$, does
the law of $\Pi u(t)$ have a $\CC^\infty$ density with respect to Lebesgue measure on $\R^N$?
\end{question}

In the finite-dimensional case, Question~\ref{mainquest} was answered successfully by H\"ormander's celebrated
`sums of squares' theorem:\cite{H1,Ho}

\begin{theorem}
Consider a collection $\{V_i\}_{i=0}^d$ of smooth vector fields on an $n$-dimensional compact manifold $\CM$
and consider the Stratonovich SDE
\begin{equation}\label{e:Strat}
dx = V_0(x)\,dt + \sum_{i=1}^dV_i(x)\circ dW_i(t)\;.
\end{equation}
Then, if the Lie algebra generated by $\d_t + V_0$ and $\{V_i\}_{i=1}^d$ spans the tangent space of $\R \times \CM$ at every point,
the law of the solutions to \eref{e:Strat} has a $\CC^\infty$ density with respect to the volume measure on $\CM$.
\end{theorem}

It turns out that a similar result still holds in the context of \eref{e:main}. Of course, the first question that needs consideration
is that of the definition of  the `Lie brackets' between the `drift vector field' $-A + F$ and the `diffusion vector fields' $q_i$.
If $G_1$ is a symmetric $k$-multilinear map and $G_2$ is a symmetric $\ell$-multilinear map, we can define a 
symmetric $k+\ell-1$-multilinear map $[G_1,G_2]$ by
\begin{equation}
[G_1,G_2](u) = \ell G_2(u,\ldots,u,G_1(u)) - k G_1(u,\ldots,u,G_2(u))\;,
\end{equation}
for every $u$ such that the right hand side makes sense. (The case $k=0$ where $G_1$ is equal to some 
constant element in $\CH$ is included.)  The problem of course is that since $A$ and $F$ are unbounded
operators in general, the domain of definition of iterated Lie brackets can rapidly shrink. We therefore need to introduce
a notion of an `admissible' Lie bracket as being one which is still defined on $\CD (A^\beta)$ for a sufficiently
low value of $\beta$. The precise definition depends on the details of the equation and can be found in \cite{Hypo}.

Note that these definitions do indeed boil down to the usual Lie brackets between vector fields with polynomial coefficients
when $\CH$ is finite-dimensional. Note also that the $k$-th iterated Lie bracket of a $k$-multilinear map with constant
elements is again a constant element. In particular, if we look at the iterated Lie brackets between
$-A + F$ and the $q_i$, many will consist of constant elements of $\CH$.
With these notions in place, our main result can be formulated as:

\begin{theorem}\label{theo:main}
In the context of \eref{e:main}, consider the collection $\Lambda \subset \CH$ of all constant elements
among the admissible iterated Lie brackets between the drift $-A + F$ and the constant elements $q_i$. (The drift is allowed to appear multiple times, so that $\Lambda$ is countably infinite in general.) 
Then, if $F$ is sufficiently regular and the solutions to \eref{e:main} are well-behaved,
Question~\ref{mainquest} has a positive answer, provided that the linear span of $\Lambda$ is dense in $\CH$. 
\end{theorem}

\begin{remark}
It seems that the first result on the regularity of finite-dimensional projections for infinite-dimensional
systems under H\"ormander-type assumptions was obtained in the linear case by Ocone\cite{Oco88:288}. These techniques
were later extended by Baudoin and Teichmann\cite{BauTei05:1765}. However, both of these works required the equation to generate a flow,
rather than a semiflow, thus excluding parabolic PDEs. 
The first regularity result of this type for parabolic SPDEs was obtained by Mattingly and Pardoux\cite{MatPar06:1742} for the 
particular setting of the 2D stochastic Navier-Stokes equations. Finally, the existence of densities (but not their regularity) was 
first obtained in a setting similar to ours by Bakhtin and Mattingly\cite{BakMat}.
\end{remark}

\begin{remark}
One can deal with the case where the non-constant Lie brackets are also included in $\Lambda$. In this case,
we obtain a collection $\Lambda(u)$ for every $u \in \CH$ and we need some weak form of uniformity of this
density. However, in most interesting examples of the form (\ref{e:main}), it is known\cite{Romito,MR2002g:76075} that considering 
the non-constant Lie brackets does not provide more information.
\end{remark}

The remainder of this note is devoted to a short explanation of some aspects of the proof of Theorem~\ref{theo:main}.

\section{Some aspects of the proof}

A natural line of attack is to try to mimic the probabilistic proof of H\"ormaner's theorem, as
obtained by Malliavin\cite{Malliavin}, Kusuoka and Stroock\cite{KSAMI,KSAMII}, Bismut\cite{Bismut}, etc.
It is a fact that, under some growth restrictions, 
if a smooth function $\Phi \colon \R^M \to \R^N$ with $M \ge N$ is such that its derivative $D\Phi(x)$ 
is of maximal rank at every point, then the image of the standard Gaussian measure under $g$ has a smooth
density with respect to Lebesgue measure. The main insight of the probabilistic proof of H\"ormander's theorem is that
this fact still holds even if $N$ is infinite, so that we can take for $\Phi$ the solution to \eref{e:Strat} viewed as
a map from Wiener space (an infinite-dimensional Gaussian space) to $\R^n$.

The problem then reduces to obtaining a moment bound for the inverse of the `Malliavin matrix' $\CM_t = D\Phi D\Phi^*$.
Let us first recall in a nutshell how this is achieved in the finite-dimensional case.
In the case of \eref{e:Strat}, the Malliavin matrix is given by
\begin{equation}\label{e:Mall}
\scal{\xi,\CM_t\xi} = \sum_{i=1}^d\int_0^t \scal{\xi,J_{s,t} V_i(x(s))}^2\,ds\;,
\end{equation}
where $J_{s,t}$ denotes the derivative of the solution map at time $t$ with respect to a change in its initial condition
at time $s < t$. In the finite-dimensional case, a key step to the analysis is to use the fact that $J_{s,t} = J_{0,t} J_{0,s}^{-1}$,
so that the invertibility of $\CM_t$ is equivalent to the invertibility of the `reduced Malliavin matrix'
\begin{equation}
\scal{\xi,\CC_t\xi} = \sum_{i=1}^d\int_0^t \scal{\xi,J_{0,s}^{-1} V_i(x(s))}^2\,ds\;,
\end{equation}
which in turn is equivalent to showing that $\P(\scal{\xi,\CC_t\xi} \le \eps) = \CO(\eps^p)$ for every $p>0$.
The point is that the integrand in the definition of $\CC_t$ is a semimartingale adapted to the filtration generated by the 
Wiener process driving the equation. The idea is then to repeatedly apply It\^o's formula to the integrand, using the fact 
that if we set $Z_V(t) = \scal{\xi,J_{0,s}^{-1} V(x(t))}$ for any smooth vector field $V$, one has
\begin{equation}
dZ_V(t) = Z_{[V_0,V]}(t)\,dt + \sum_{i=1}^d Z_{[V_i,V]}(t)\circ dW(t)\;.
\end{equation}
Norris' lemma\cite{Norr} (a quantitative version of the Doob-Meyer decomposition theorem) then ensures that
if $Z_V$ is small for some $V$, it must also be small for $\{[V,V_j]\}_{j=0}^d$. Since H\"ormander's condition
ensures that these quantities cannot all be small simultaneously, the claim follows.

The problem with the generalisation to parabolic SPDEs comes from the fact that $J_{0,s}$ is not invertible, so that 
we have to deal with the non-adapted integrand in \eref{e:Mall}. The trick is to use the fact that all the vector fields
appearing in our case have a polynomial structure. This allows us to exploit the additive structure of our noise by 
setting $v = u - QW$, so that $v$ satisfies the random PDE
\begin{equation}
{dv \over dt} = -Av + F(v+ QW) + AQW\;.
\end{equation}
Since $F$ is polynomial, the right hand side of this equation can be written as a sum over finitely many terms,
each of them being multilinear in both $v$ and $W$. A similar procedure to the finite-dimensional analysis then
allows us to reduce the question to the following version of Norris' lemma. 
For a multiindex $\alpha = (\alpha_1,\ldots,\alpha_\ell)$ write
$W_\alpha(t) = W_{\alpha_1}(t)\cdot\ldots\cdot W_{\alpha_\ell}(t)$ (with the convention that $W_\emptyset(t) = 1$)
and consider a stochastic process of the form
\begin{equation}\label{e:Z}
Z(t) = \sum_{|\alpha| \le m} A_\alpha(t) W_\alpha(t)\;,
\end{equation}
where the $A_\alpha$ are stochastic processes that are not necessarily adapted to the Brownian filtration, but that
are almost surely Lipschitz continuous in time. Then, one has:

\begin{proposition}\label{prop:Wiener}
In the above setting, there exists a universal family of events $\Omega_\eps$ depending only on $m$
such that $\P(\Omega_\eps) = \CO(\eps^p)$ for every $p$ and such that the implication
  \begin{equation}\label{e:dichotomy}
    \|Z\|_{L^\infty} \le \eps \quad\Longrightarrow\quad
    \left\{\begin{array}{rl} \text{either}& \sup_\alpha
        \|A_\alpha\|_{L^\infty} \le \eps^{3^{-m}} \\ \text{or} &
        \sup_\alpha \|A_\alpha\|_{\Lip} \ge
        \eps^{-3^{-(m+1)}} \end{array}\right.
  \end{equation}
 holds for every $W \not\in \Omega_\eps$ and for every $\eps \in (0,1]$. 
 Here, the supremum norms are taken over the time 
 interval $[0,1]$.
\end{proposition}

\begin{remark}
Note that for any given $W \not \in \Omega_\eps$, \eref{e:dichotomy} is a \textit{deterministic}
implication that holds simultaneously for \textit{all} processes of the form (\ref{e:Z}).
\end{remark}

In order to be able to make use of this proposition, it is important to note that if $G$ is any polynomial 
 map (\i.e. a sum of finitely many multilinear maps),
then there exist $m>0$ and finitely many polynomial maps $G_\alpha$ that all consist of iterated Lie brackets between
$G$, $-A+F$, and the $q_i$ such that the identity
\begin{equation}
{d\over ds}\scal{\xi,J_{s,t} G(v(s))} = \sum_\alpha \scal{\xi,J_{s,t} G_\alpha(v(s))} W_\alpha(t)
\end{equation}
holds. Using Proposition~\ref{prop:Wiener} repeatedly, we conclude that if $\scal{\xi,\CM_t\xi}$ is small, then
$\scal{\xi,J_{s,t} G(v(s))}$ must be small for every admissible Lie bracket $G$ constructed from the drift and the
diffusion coefficients of our original SPDE. Since these span a dense linear subspace of $\CH$, this cannot be
true, so that any finite-dimensional projection of $\CM_t$ must indeed be invertible, thus concluding the proof.

\bibliographystyle{./ws-procs9x6}
\bibliography{./refs}
\end{document}